\documentclass[12pt, a4paper, preprintnumbers, showkeys]{revtex4} 

\usepackage{tikz}
\usepackage[fleqn]{amsmath}  
\usepackage{amssymb}   
\usepackage{mathrsfs} 
\usepackage{bm}        
\usepackage{graphicx}  
\usepackage{subfigure} 
\usepackage{pict2e}   
\usepackage{ulem}     
\usepackage{cases}
\usepackage[pdfstartview=FitH,
            bookmarksnumbered=true,
            bookmarksopen=true,
            colorlinks=true,  
            pdfborder=001,    
            linkcolor=blue,   
            citecolor=blue,   
            urlcolor=blue    
            ]{hyperref}      
\usepackage[left=2.5cm,  right=2.5cm,  top=2.5cm,  bottom=2.2cm]{geometry}
\usepackage{ntheorem}

\allowdisplaybreaks[4]
\newcommand*{\QED}{\hfill\ensuremath{\square}}  

\newtheorem{lemma}{\hskip 24pt Lemma}[section] 
\newtheorem{theorem}{\hskip 24pt Theorem}[section]

\newtheorem{definition}{\hskip 24pt Definition}[section]
\newtheorem{claim}{\hskip 24pt Claim}[section]

\begin{document}


 \title{Extremal spectral results of planar graphs without
 $C_{l, l}$ or $\mathrm{Theta}$ graph}

\author{Hao-Ran Zhang, Wen-Huan Wang\footnote{Corresponding author. Email: whwang@shu.edu.cn}}

\affiliation{Department of Mathematics,  Shanghai University,  Shanghai 200444,  China}
\affiliation{Newtouch Center for Mathematics of Shanghai University, Shanghai 200444, China}
\date{\today}

\begin{abstract}
	
   Let  $\mathcal {F}$ be a given family of graphs.
     A graph  $G$ is $\mathcal {F}$-free if it does not contain any member of $\mathcal {F}$  as a subgraph.
  Let $C_{l, l}$ be a graph obtained from $2C_l$ such that the two cycles share a common vertex, where  $l\geqslant3 $.
   A $\mathrm{Theta}$ graph is obtained from a cycle $C_k$ by adding an additional edge between two non-consecutive vertices on $C_k$, where  $k\geqslant 4$.
     Let $\Theta_k$ be the set of $\mathrm{Theta}$ graphs on $k$ vertices, where  $k \geqslant 4$.
      For sufficiently large $n $, the unique extremal planar graph with the maximum spectral radius among  $C_{l, l}$-free planar graphs on $n$ vertices and
       among $\Theta_k$-free  planar graphs on $n$ vertices are characterized respectively,  where $l \geqslant 3$ and $k \geqslant 4$.

   \end{abstract}


\maketitle


\section{\label{introduction}Introduction}

     Let $G=(V(G), E(G))$ be a graph, where $V(G)$ and $E(G)$ are  the sets of vertices and edges of $G$, respectively.
      Let  $\mathcal {F}$ be a given family of graphs.
     If a graph $G$ does not contain any member of $\mathcal {F}$  as a subgraph, then $G$ is $\mathcal {F}$-free.
       If $\mathcal {F}=\{F\}$, then $G$ is $F$-free.
       The Tur\'{a}n number  of $\mathcal {F}$, denoted by $\textrm{ex}(n, \mathcal{F})$, is the maximum number of the edges in  $\mathcal {F}$-free graphs on $n$ vertices.
       The  Tur\'{a}n type extremal problem is  to determine $\textrm{ex}(n, \mathcal{F})$.
      Studying the Tur\'{a}n number $\textrm{ex}(n, \mathcal{F})$  is a central
    problem in extremal graph theory.
     In 2010, Nikiforov   \cite{NIKIFOROV20102243} proposed the following  spectral analogue of the Tur\'{a}n type problem for graphs:
   ``What is the maximum spectral radius of graphs of  order $n$, not containing a given graph $F$?"
    The study of the spectral Tur\'{a}n type extremal problem have drawn increasingly extensive interest from researchers.
         For  the results about the Tur\'{a}n number and the spectral Tur\'{a}n type problem for graphs, one can refer to surverys \cite{Furedi1991,Mubayi2016}          and
          Refs.\ \cite{GAO2019, ZHAI2021103322, CIOABA2022103420,  2023-JingWang-p20}.

    Let $P_k=v_1v_2\cdots v_k$ and $C_k=v_1v_2\cdots v_kv_1$ be a  path and a cycle with $k$ vertices, respectively,  where $k\geqslant 3$.
    Let $tC_l$ be the disjoint union of $t$ copies of $C_l$, where $t\geqslant 1$ and $l\geqslant3 $.
    A $\mathrm{Theta}$ graph  is obtained from a cycle $C_k$ by adding an additional edge between two non-consecutive vertices on $C_k$, where  $k\geqslant 4$.
    Let $C_a\cdot C_b$ be the graph obtained from $C_a$ and $C_{b}$ by identifying $u_1$ with $v_1$ and identifying $u_2$ with $v_2$, where $u_1u_2\in E(C_a)$, $v_1v_2\in E(C_b)$ and $a,  b\geqslant 3$.
    Obviously, $C_a\cdot C_b$ is a $\mathrm{Theta}$ graph.
     Let  $\Theta_k$ be the set of all the $\mathrm{Theta}$ graphs on $k$ vertices,  where $k\geqslant4 $.
    Namely $\Theta_k=\{C_3\cdot C_{k-1}, C_4\cdot C_{k-2}, \cdots, C_{\lfloor\frac{k}{2}\rfloor+1}\cdot C_{\lceil\frac{k}{2}\rceil+1}\}$,  where $k\geqslant4 $.
   When $k=4$, we have $\Theta_4=\{C_3\cdot C_{3}\}$.
    For two given graphs $G_1$ and $G_2$,  the joint product of $G_1$ and $G_2$, denoted by  $G_1+G_2$, is  obtained from $G_1$ and $G_2$ by adding an edge between an arbitrary vertex in $G_1$ and an arbitrary vertex in $G_2$.

     More recently, extremal problems have been considered where the host graph is taken
   to be a planar graph.
    For  a given planar graph family $\mathcal {F}$,  the planar  Tur\'{a}n number,  denoted by $\textrm{ex}_ {P}(n, \mathcal {F})$,  is the maximum number of the edges over all
    $\mathcal {F}$-free planar graphs on $n$ vertices. In 2016, Dowden  \cite{2016-Dowden-p213} initiated  the study of  $\textrm{ex}_ {P}(n, \mathcal {F})$   and
obtained  tight bounds for $\textrm{ex}_{P}(n, C_4)$ and $\textrm{ex}_{P}(n, C_5)$.
  Lan et al.\ \cite{2019-Lan-111610} and Ghosh et al.\ \cite{doi:10.1137/21M140657X} characterized the bounds for $\textrm{ex}_ {P}(n, C_6)$.
 Lan et al.\ \cite{2019-Lan-111610} determined the bounds for $\textrm{ex}_{P}(n,\Theta_k)$ with $k=4, 5, 6$.
  Fang et al.\  \cite{FANG2022112794} obtained  $\textrm{ex}_ {P}(n,H_k)$ and $\textrm{ex}_{P}(n,F_k)$, 
   where $H_k$ is a  friendship graph with $2k+1$ vertices which is obtained from $kC_3$ such that all the $C_3$ share a common vertex and  $F_k$ is a $k$-fan with $F_k\cong K_1+P_{k+1}$.

    Studying  the spectral radius of planar graphs has a long history.
      Boots and Royle \cite{articleboots} and  independently Cao and Vince \cite{CAO1993251} gave such a conjecture: $P_2+P_{n-2}$ is the unique planar graph with the maximum spectral radius, where $n\geqslant 9 $.
     In 2017, Tait and Tobin \cite{TAIT2017137} proved that this conjecture is right when $n$ is sufficiently large.

      The planar extremal spectral radius,  denoted by $\textrm{spex}_ {P}(n, \mathcal {F})$, is the maximum spectral radius over all $\mathcal {F}$-free planar graphs  on $n$ vertices.
    An $\mathcal {F}$-free planar graph on $n$ vertices with maximum spectral radius is called an extremal graph  to $\textrm{spex}_{P}(n, \mathcal {F})$.
   Recently, researchers start to study $\textrm{spex}_{P}(n, \mathcal {F})$ and some results are obtained.
     In 2022, Zhai and Liu \cite{mingqing2022extremal} first determined the maximum number of edges in planar graphs of order $n$ and maximum degree $n-1$ without $k$ edge-disjoint cycles, and then obtained the maximum spectral radius as well as its unique extremal graph over all planar graphs on $n$ vertices without $k$ edge-disjoint cycles.
   In 2024, Fang et al.\ \cite{https://doi.org/10.1002/jgt.23084} characterized the unique extremal planar graph to  $\textrm{spex}_{P}(n, tC_l)$ and    $\textrm{spex}_{P}(n, t\mathscr{C})$,   where $t\mathscr{C}$ are the family of $t$ vertex-disjoint  cycles without length restriction, $t\geqslant 1$,  $l\geqslant 3$, and  $n$ is sufficiently large.

    Let $C_{l, l}$ be a graph obtained from $2C_l$ such that the two cycles share a common vertex, where  $l\geqslant3 $.
    Inspired by the results obtained  as above,   in this article, we will study the  planar extremal spectral radius in  planar graphs on  $n$ vertices  without $C_{l, l}$ or without any member of $\Theta_k$,
     where $l \geqslant 3$  and $k \geqslant 4$.
     When $n $ is sufficiently large,  the unique extremal graph  to  $\textrm{spex}_{P}(n, C_{l, l})$  and  $\textrm{spex}_{P}(n, \Theta_k)$  are characterized respectively, where $l \geqslant 3$  and $k \geqslant 4$.
   The results are shown in Theorems \ref{MR1}--\ref{MR3aa}.

 \section{\label{Preliminaries}Preliminaries}

  To deduce the main results of the present paper, some definitions and necessary lemmas  are simply quoted here.

    Let $G=(V(G), E(G)) $  be a graph.
   For a vertex $v\in V (G)$, let $N(v)$ be the set of the vertices which are adjacent to $v$.
   The degree of $v$, denoted by  $d_G(v)$,  is the number of the edges which are incident with $v$. For convenience, we denote $d_G(v)$ by $d(v)$.

    Let $A(G)$ be the adjacency matrix of $G$.
     The spectral radius of $G$, denoted by $\rho(G)$, is the largest modulus of all the eigenvalues of $A(G)$.
      Let $\bm{x} =(x_{1},x_{2},\cdots,x_{n})^{\mathrm{T}}$ be  an eigenvector of $A(G)$ corresponding to $\rho(G)$.
        The eigenequation  for any vertex $v$ of $G$ is
 \begin{align}\label{q1}
	\rho(G)x_v=\sum \limits_{uv\in E(G)}x_u.
  \end{align}
   The Rayleigh quotient of  $G$ is expressed as:
  \begin{align}\label{q2}
	\rho(G)= \max \limits_{\bm{x}\in \mathbb{R}^n}\frac{\bm{x}^\mathrm{T}A(G)\bm{x}}{\bm{x}^\mathrm{T}\bm{x}}
	=\max \limits_{\bm{x}\in \mathbb{R}^n}\frac{2\sum\limits_{uv\in E(G)}x_ux_v}{\bm{x}^\mathrm{T}\bm{x}},
 \end{align} where  $\mathbb{R}^n$ is the set of  $n$-dimensional real numbers.

 If $G$ is a connected graph with $n$ vertices, then by the well-known Perron-Frobenius theorem, there exists a unique positive vector $\bm{x}=(x_1, \cdots, x_n)^\mathrm{T}$ corresponding  to $\rho(G)$.
      It is  convenient for us to normalize $\bm{x}$ such that the maximum entry of $\bm{x}$ is 1.
      Usually,  $\bm{x}$ is referred to  as the Perron vector of $G$.

  \begin{lemma}\label{PF1} \cite{https://doi.org/10.1002/jgt.23084}
	Let $F$ be a  planar graph and $n \geqslant \max\{2.16\times 10^{17}, 2|V(F)|\}$.
    If $F$ is a subgraph of $2K_1+P_{n/2}$, but not of $K_{2, n-2}$,  then the extremal graph having the maximum spectral radius among $F$-free planar graphs  contains a copy of $K_{2, n-2}$.
	\end{lemma}


%
%
%
 \begin{definition}\label{df1}\cite{https://doi.org/10.1002/jgt.23084}
	Let $s_1$ and $s_2$ be two positive integers, where $s_1\geqslant s_2 \geqslant 1$.
 Let $H=P_{s_1}\cup P_{s_2}\cup H_0$, where $H_0$ is a union of vertex-disjoint paths.
 If 	$$
		H^*=\begin{cases}
			P_{s_1+1}\cup P_{s_2-1}\cup H_0, &s_2\geqslant 2, \\
			P_{s_1+s_2}\cup H_0, &s_2=1,
		\end{cases}
		$$
	then we say that $H^*$ is a graph obtained from $H$ by  a $(s_1, s_2)$-transformation.
 \end{definition}

  Obviously, $H^*$ is a union of vertex-disjoint paths too and $K_2+H^*$ is a planar graph.


 \begin{lemma}\cite{Li1979}
 \label{liqiao}
 Let $G$ be a connected graph, and let $G'$ be a proper spanning subgraph of $G$.
 Then  $\rho(G)> \rho(G')$. \end{lemma}

 \section{Main results }

 \subsection{Extremal spectral results of planar graphs without $C_{l, l}$}

  In this section,  among the set of planar graphs on $n$ vertices without $C_{l, l}$,  we characterize the graph with the maximum spectral radius,  where $l\geqslant3 $ and $n $ is sufficiently large.
  The results are shown in Theorems \ref{MR1} and \ref{MR2}.

  Let $n$,  $n_1$ and  $n_2$ be three positive integers,  where $n > n_1>n_2$.
  We define a graph $H(n_1, n_2)$ as follows. Let
  $$
  H(n_1, n_2)=\begin{cases}
	P_{n_1}\cup \frac{n-2-n_1}{n_2}P_{n_2}, & if~ n_2 | (n-2-n_1), \\
	P_{n_1}\cup \lfloor \frac{n-2-n_1}{n_2} \rfloor P_{n_2}\cup P_{n-2-n_1-\lfloor \frac{n-2-n_1}{n_2} \rfloor n_2}, &\text{otherwise}.
 \end{cases}
 $$
 Specifically,  when $n_1=n_2=1$,  $H(1, 1)$ is a graph having  $n-2$ isolated vertices.

  For a graph $G$ and  a  vertex subset $W \subseteq V(G)$, let $G[W]$ be the subgraph induced by $W$ whose
edges are precisely the edges of $G$ with both ends in $W$.

 \begin{theorem}\label{MR1}
   Let $G$ be a planar graph on $n$ vertices without $C_{3, 3}$.
  When $n\geqslant 2.16\times 10^{17}$,  we  have
   $\rho(G)\leqslant \rho(K_2+H(1, 1))$  with the equality if and only if $G\cong K_2+H(1, 1)$.
\end{theorem}

\textbf{Proof.}
  We suppose that $G^*$ is the  extremal graph  to $\textrm{spex}_{P}(n, C_{3, 3})$,  where $n\geqslant 2.16\times 10^{17}$.
  Since $C_{3, 3}\subseteq2K_1+P_{n/2}$ and $C_{3, 3} \nsubseteq K_{2, n-2}$,  by Lemma \ref{PF1},   we have $ K_{2, n-2}\subseteq G^*$.
  By Perron-Frobenius theorem,  there exists a positive eigenvector $\bm{x}=(x_1,  x_2,  \cdots,  x_n)^{\mathrm{T}}$ corresponding to $\rho(G^*)$.
    By the proof of Lemma \ref{PF1} in Ref.\  \cite{https://doi.org/10.1002/jgt.23084} (see Pages 6-7),  among $\bm{x}$,  we know that there exist two components with larger values.
   In $G^*$,  let $u'$ and $u''$ be the two vertices corresponding to the two larger components,  where $u'$ satisfies $x_{u'}=\max\{x_i: 1\leqslant i \leqslant n\}=1$.
   Let $D=\{u', u''\}$.
  In $G^*$,  let $R=V(G^*)\setminus D=\{v_1, \cdots , v_{n-2}\}$.
  Since $K_{2, n-2}\subseteq G^*$, we have $V(G^*)=V(K_{2, n-2})$ and
  $N(u')\cap N(u'')=R$ (by the proofs of Lemma \ref{PF1} in Ref.\ \cite{https://doi.org/10.1002/jgt.23084}).

  We prove $|E(G^*[R])\cup E(G^*[D])|\leqslant 1$.
  Otherwise,  we suppose  $|E(G^*[R])\cup E(G^*[D])|\geqslant 2$.
  Let $e_1, e_2 \in E( G^*[R]) \cup E(G^*[D])$. Two cases are considered as follows.

  \textbf{Case (1)}. $e_1=u'u''$,  $e_2\in E(G^*[R])$.

   Without loss of generality,  we suppose $e_2=v_1v_2$.
   In $G^*$,  there exists a subgraph $C_{3, 3}$ formed by two cycles $u'v_1v_2u'$ and $u'u''v_3u'$.

  \textbf{Case (2)}. $e_1, e_2\in E(G^*[R])$.

  If  $e_1\cap e_2=\varnothing$,  we suppose $e_1=v_1v_2$ and  $e_2=v_3v_4$. We can find a subgraph $C_{3, 3}$  in $G^*$ formed by two cycles $u'v_1v_2u'$ and $u'v_3v_4u'$.
   If   $e_1\cap e_2\neq\varnothing$,  we suppose $e_1=v_1v_2$ and  $e_2=v_2v_3$. Obviously,   $u'v_1v_2u'$ and $u''v_2v_3u''$ form a subgraph $C_{3, 3}$ in $G^*$.

   By combining the proofs of Cases (1) and (2),  we get that $G^*$ contains a subgraph $C_{3, 3}$.  This is a contradiction.
   Thus,  we get $|E(G^*[R])\cup E(G^*[D])|\leqslant 1$.

   Next,  we prove $G^*\cong K_2+H(1, 1)$.
   Since  $|E(G^*[R])\cup E(G^*[D])|\leqslant 1$,  we have  $G^*\cong K_2+H(1, 1)$ or $G^*\cong K_{2, n-2}$ or $G^*\cong K^+_{2, n-2}$,
   where $K^+_{2, n-2}=K_{2, n-2}+v_iv_j$ with $v_i,  v_j\in V(K_{2, n-2})\setminus \{u', u''\}=\{v_1,  \cdots,  v_{n-2}\}$ and  $1\leqslant i<j\leqslant n-2$.
   We suppose  $v_i v_j=v_1v_2$.
    Since $K_{2, n-2}\subset K_2+H(1, 1)$,  by Lemma \ref{liqiao},  $\rho(K_{2, n-2})< \rho(K_2+H(1, 1))$.
   Thus, $G^*\cong K_2+H(1, 1)$ or $G^*\cong K^+_{2, n-2}$.
    We suppose $G^*\cong K^+_{2, n-2}$.
    Among $K^+_{2, n-2}$,  by the symmetry of $u'$ and  $u''$, $x_{u''}=x_{u'}=1$.
    	 Among $K^+_{2, n-2}$, for any $v_i\in R$, we have $x_{v_i} \leqslant \frac{1}{10}$ (by Claim 2.6 on Page 7 in Ref.\ \cite{https://doi.org/10.1002/jgt.23084}), where $1\leqslant i \leqslant n-2$.
 Namely $x_{v_1}, x_{v_2}\leqslant \frac{1}{10}$.
    Therefore, we obtain
     \begin{align} \label{q1jiwwa}
     \rho (K^+_{2, n-2})-\rho (K_2+H(1, 1))&\leqslant \frac{\bm{x}^\mathrm{T}\left(A(K^+_{2, n-2})-A(K_2+H(1, 1))\right)\bm{x}}{\bm{x}^\mathrm{T}\bm{x}}  \nonumber\\
     &=\frac{2}{\bm{x}^\mathrm{T}\bm{x}}(x_{v_1}x_{v_2}-x_{u'}x_{u''})\leqslant \frac{2}{\bm{x}^\mathrm{T}\bm{x}}(\frac{1}{100}-1)<0.
      \end{align}
    Namely, $\rho (K^+_{2, n-2})<\rho (K_2+H(1, 1))$.
    This contradicts $G^*\cong K^+_{2, n-2}$.
    Thus, $G^*\cong K_2+H(1, 1)$.\QED

\begin{theorem}\label{MR2}
	 Let $G$ be a planar graph on $n$ vertices without $C_{l, l}$,  where $l\geqslant 4$.
  When $n\geqslant \max\{2.16\times 10^{17}, 9\times 2^{l-1}+3, \frac{625}{32}\lfloor\frac{l-3}{2}\rfloor^2+4\}$, we have  $\rho(G)\leqslant \rho(K_2+H(l-2, l-2))$
  with the equality if and only if $G\cong K_2+H(l-2, l-2)$.
\end{theorem}

 \textbf{Proof.}
   We suppose that $\widetilde{G}$ is  the  extremal graph  to $\textrm{spex}_{P}(n, C_{l, l})$,  where $l\geqslant 4$ and  $n\geqslant \max\{2.16\times 10^{17}, 9\times 2^{l-1}+3, \frac{625}{32}\lfloor\frac{l-3}{2}\rfloor^2+4\}$.
   Since $n\geqslant 9\times 2^{l-1}+3> 2(2l-1)=2|V(C_{l, l})|$, we  get $C_{l, l}\subseteq 2K_1+P_{n/2}$ and $C_{l, l} \nsubseteq K_{2, n-2}$.
    By Lemma \ref{PF1}, $K_{2, n-2}\subseteq \widetilde{G}$.
   Thus,  by Lemma \ref{liqiao}, $\rho(\widetilde{G})\geqslant \rho(K_{2, n-2})=\sqrt{2n-4}$.
   By Perron-Frobenius theorem,  we get that there exists a positive eigenvector $\bm{x}=(x_1,  x_2,  \cdots,  x_n)^{\mathrm{T}}$ corresponding to $\rho(\widetilde{G})$.
     By the proofs of Lemma \ref{PF1} in Ref.\  \cite{https://doi.org/10.1002/jgt.23084} (see Pages 6-7),  among $\bm{x}$,  we know that there exist two components with larger values.
   In $\widetilde{G}$,  let $u'$ and $u''$ be the two vertices corresponding to the two larger components,  where $u'$ satisfies $x_{u'}=\max\{x_i: 1\leqslant i \leqslant n\}=1$. Let $D=\{u', u''\}$.
   In $\widetilde{G}$,  let $R=V(\widetilde{G})\setminus D=\{v_1, \cdots , v_{n-2}\}$.
     Since   $K_{2, n-2}\subseteq \widetilde{G}$,  we get $R=N(u')\cap N(u'')=V(\widetilde{G})\setminus D$  (by the proofs of Lemma \ref{PF1} in Ref.\ \cite{https://doi.org/10.1002/jgt.23084}).

    To obtain Theorem \ref{MR2},  Claims \ref{Cl1}--\ref{Cl77} are  needed.

 \begin{claim}\label{Cl1}
	$ \widetilde{G}[R]$ is a union of $t$ vertex-disjoint paths, where $t\geqslant 2$.
	\end{claim}

  \textbf{Proof.}
  To obtain Claim \ref{Cl1}, we only need to prove that $ \widetilde{G}[R]$ does not contain cycles and the maximum degree of each vertex  of  $ \widetilde{G}[R]$ is less than or equal to $2$.
  We assume that $ \widetilde{G}[R]$ contains a cycle $C_k$, where $3\leqslant k\leqslant n-2$.
   Let $C_k=v_1v_2\cdots v_kv_1$.
   If $3\leqslant k<n-2$, then there exists an $i\in \{k+1, \cdots , n-2\}$ such that $v_i\notin V(C_k)$.
   We suppose $i=n-2$.
   Then $\widetilde{G}[V(C_k)\cup \{u', u'', v_{n-2}\}]$ contains  $K_5$-minor as a subgraph.
   This contradicts the fact that $\widetilde{G}$ is a planar graph.
    If $k=n-2$, then $\widetilde{G}$ contains  $C_{l, l}$  as a subgraph.
    This is a contradiction. Thus, we get that $ \widetilde{G}[R]$ does not  contain cycles.
   We suppose that there exists a vertex (denoted by $v_1$) in $ \widetilde{G}[R]$ such that $d_{ \widetilde{G}[R]}(v_1)\geqslant 3$.
   Therefore, there exist  at least three edges in $ \widetilde{G}[R]$  such that they are  incident with $v_1$.
   We suppose that  $v_1v_2, v_1v_3, v_1v_4\in E( \widetilde{G}[R])$.
   Then $\widetilde{G}[\{u', u'', v_1, v_2, v_3, v_4\}]$ contains $K_{3, 3}$.
   This contradicts the fact that $\widetilde{G}$ is a planar graph.
  Therefore,  the maximum degree of each vertex  of  $ \widetilde{G}[R]$ is less than or equal to $2$.
   We conclude that $ \widetilde{G}[R]$ is a union of $t$ vertex--disjoint paths, where $t\geqslant 1$.
   If $t=1$, then $ \widetilde{G}[R]=P_{n-2}$.
    Since $n\geqslant \max\{2.16\times 10^{17}, 9\times 2^{l-1}+3, \frac{625}{32}\lfloor\frac{l-3}{2}\rfloor^2+4\}$,
   obviously, $\widetilde{G}$ contains  $C_{l, l}$ as a subgraph. This is a contradiction. Thus, $t\geqslant 2$.\QED

  \begin{claim}\label{Cl2}
	$u'u''\in E(\widetilde{G})$.
\end{claim}

\textbf{Proof.}
    Since $ \widetilde{G}[R]$ is a union of $t$ vertex--disjoint paths (by Claim \ref{Cl1}),
    there exists at least an $i\in\{1, \cdots, n-2\}$  such that $v_iv_{i+1}\notin E( \widetilde{G}[R])$ (we define $v_{n-1}=v_1$), where $t\geqslant 2$.
     Thus, $u'v_iu''v_{i+1}u'$  is a face of $\widetilde{G}$, where $1 \leqslant i \leqslant n-2$.
     We suppose $u'u''\notin E(\widetilde{G})$.
    Let $G^{'}=\widetilde{G}+\{u'u''\}$, where $u'u''$ passes through the surface $u'v_iu''v_{i+1}u'$ with $1 \leqslant i \leqslant n-2$.
    Obviously, $G^{'}$ is  a planar graph.
    We prove that $G^{'}$ does not contain  $C_{l, l}$ as a subgraph.
    Otherwise, we suppose that $G^{'}$ contains  $C_{l, l}$ as a subgraph.
    Let $C^1$ and $C^2$  be the two cycles of $C_{l, l}$.
    Obviously, the length of $C^1$ and $C^2$ is $l$.
    Since $ \widetilde{G}$ does not contain $C_{l, l}$, we get $u'u''\in E(C^1)\cup E(C^2)$.
    Suppose $u'u''\in E(C^1)$.  Let $C^1=u'u''v_1v_2\cdots v_{l-2}u'$.  Three cases are considered as follows.

    \textbf{Case (1)}. $V(C^1)\cap V(C^2)=\{u'\}$.

  We suppose $C^2=u'v_{l-1}v_l\cdots v_{2l-3}u'$. Obviously, among  $\widetilde{G}$,  $u'v_1u''v_2\cdots v_{l-2}u'$ and $C^2$ are two cycles of length $l$ and they intersect at $u'$.
  Namely,  $\widetilde{G}$ contains $C_{l, l}$ as a subgraph. This is a contradiction.

 \textbf{Case (2)}. $V(C^1)\cap V(C^2)=\{u''\}$.

   By the symmetry of $u'$ and $u''$ and the proofs of Case (1),  we  can also  get a contradiction in Case (2).

 \textbf{Case (3)}. $V(C^1)\cap V(C^2)=\{v_i\}$, where $i\in \{1, 2, \cdots, l-2\}$.

  Without loss of generality, we suppose $V(C^1)\cap V(C^2)=\{v_{l-2}\}$.
  We can check that $C^2$ is contained in $ \widetilde{G}[R]$.
  We assume $C^2=v_{l-2}v_{l-1}\dots v_{2l-3}v_{l-2}$.
  It is checked that $\widetilde{G}[\{u', u''\}\cup V(C^1)\cup V(C^2)]$ contains $K_{3, 3}$-minor as a subgraph.
  This contradicts the fact that $\widetilde{G}$ is a planar graph.

  By combining the proofs of Cases (1)--(3),  we get that $G^{'}$ is a planar graph without $C_{l, l}$, where  $l\geqslant 4$.
  Since  $\widetilde{G}\subset G^{'}=\widetilde{G}+\{u'u''\}$, by Lemma  \ref{liqiao}, we have $\rho(\widetilde{G}) < \rho(G^{'})$. This contradicts the definition of $\widetilde{G}$.
  Thus, $u'u''\in E(\widetilde{G})$.\QED

\begin{claim}\label{Cl3}
	For any $v\in R$, we have $x_v\in[\frac{2}{\rho(\widetilde{G})}, \frac{2}{\rho(\widetilde{G})}+\frac{6}{\rho^2(\widetilde{G})}]$.
\end{claim}

\textbf{Proof.}  By Claim \ref{Cl2}, among $ \widetilde{G}$, we get that  $u'$ and $u''$ are symmetric. Therefore,  we have $x_{u'}=x_{u''}=1$.
 For any vertex  $v\in R$, using the eigenequation (\ref{q1}) for $v$ of $\widetilde{G}$, we get
 \begin{align}\label{q1jia}
 	\rho(\widetilde{G}) x_v=x_{u'}+x_{u''}+\sum \limits_{u\in N_R(v)}x_u=2+\sum \limits_{u\in N_R(v)}x_u.
 \end{align}

  Since  $ \widetilde{G}[R]$ is a union of $t$ ($t\geqslant 2$) vertex--disjoint paths (by Claim \ref{Cl1}), for any vertex  $v\in R$, $d_R(v)\leqslant 2$.
  For any vertex $v\in R$, we have $x_{v} \leqslant \frac{1}{10}$ (by Claim 2.6 on Page 7 in Ref.\ \cite{https://doi.org/10.1002/jgt.23084}).
 Thus, it follows from (\ref{q1jia}) that $2\leqslant \rho(\widetilde{G}) x_v\leqslant 2+2\times\frac{1}{10}=2+\frac{1}{5}$.
 Namely, for any $v\in R$, $\frac{2}{\rho(\widetilde{G})}\leqslant x_v\leqslant \frac{2}{\rho(\widetilde{G})}+\frac{1}{5\rho(\widetilde{G})}\leqslant \frac{3}{\rho(\widetilde{G})}$.
  Using (\ref{q1jia}) again,  we get $2\leqslant \rho(\widetilde{G}) x_v\leqslant 2+2\times\frac{3}{\rho(\widetilde{G})}\leqslant 2+\frac{6}{\rho(\widetilde{G})}$.
  Thus, $\frac{2}{\rho(\widetilde{G})}\leqslant x_v\leqslant \frac{2}{\rho(\widetilde{G})}+\frac{6}{\rho^2(\widetilde{G})}$.
  \QED

 \begin{claim}\label{Cl44}
	   Let $H$ and $H^*$ be the two graphs as shown in
	 Definition \ref{df1}, $n \geqslant max\{2.16\times 10^{17}, 9\times 2^{s_2+1}+3,2|V(C_{l,l})|\}$. When $ \widetilde{G}[R]\cong H$, we have $\rho(K_2+H^*)>\rho(\widetilde{G})$.
		\end{claim}
	
\textbf{Proof.} By Claims \ref{Cl1}--\ref{Cl3} and using the same methods as those for the proofs of Lemma 3.2 in  Ref.\ \cite{https://doi.org/10.1002/jgt.23084}  (see Pages 10-13), we obtain Claim \ref{Cl44}.\QED

  Let $H$ be a  union of $t$  vertex--disjoint paths, where $t\geqslant 2$.
  For any $i\in \{1, \cdots, t\}$, let $n_i(H)$ be the order of the $i$-th longest path of $H$.
  Obviously, $n_1(H)\geqslant n_2(H)\geqslant \cdots \geqslant n_t(H)$.

 \begin{claim}\label{Cl4}
	Let $H$ be a  union of $t$  vertex--disjoint paths, where $t\geqslant 2$.
 Then $K_2+H$ does not contain $C_{l, l}$  if and only if $n_1(H)+n_2(H)\leqslant 2l-4$ and $n_1(H)\leqslant l-2$ or $n_1(H)+n_2(H)\leqslant 2l-4$ and $n_2(H)+n_3(H)\leqslant l-3$.
\end{claim}

\textbf{Proof.}  To obtain  Claim \ref{Cl4},  we only need to prove that $K_2+H$ contains $C_{l, l}$ if and only if $n_1(H)+n_2(H)\geqslant 2l-3$,  or  $n_1(H)+n_2(H)\leqslant 2l-4$, $n_1(H)\geqslant l-1$ and $n_2(H)+n_3(H)\geqslant l-2$.
 Let $V(K_2)=\{u', u''\}$ and
 $V(H)=\{v_i: 1 \leqslant i \leqslant |V(H)|\}$.

   (1). The proof of sufficiency of  Claim \ref{Cl4}.

    Two cases are considered as follows.

 \textbf{Case (1)}.  $n_1(H)+n_2(H)\geqslant 2l-3$.

 Since $n_1(H)\geqslant n_2(H)$, we have $n_1(H)\geqslant l-1$.

\textbf{Subcase (1.1)}.  $n_1(H)\geqslant 2l-3$.

 When $n_1(H)\geqslant 2l-3$,  there exists a subgraph $P_{2l-3}$ in $H$.
 We suppose $P_{2l-3}=v_1v_2\cdots v_{2l-3}$.
 Then $K_2+H$ contains  $C_{l, l}$ as a subgraph, where  $C_{l, l}$ is formed by two cycles $u'v_1v_2\cdots v_{l-1}u'$ and $u'v_{l}v_{l+1}\cdots v_{2l-4}u''v_{2l-3}u'$ which are of length $l$ and they intersect at $u'$.

\textbf{Subcase (1.2)}. $l-1\leqslant n_1(H)\leqslant 2l-4$.

  First, we assume $n_1(H)+n_2(H)= 2l-3$.
  Since $n_1(H)\geqslant n_2(H)$, we get $l-1=\lceil \frac{2l-3}{2} \rceil \leqslant n_1(H)\leqslant  2l-4$
   and  $1 \leqslant n_2(H)\leqslant  \lfloor \frac{2l-3}{2} \rfloor=l-2$. 

  When $n_1(H)=l-1$ and  $n_2(H)=l-2$, in $H$, there exist two subgraphs $P_{l-1}$ and $P_{l-2}$.
  We suppose $P_{l-1}=v_1v_2\cdots v_{l-1}$ and  $P_{l-2}=v_{l}v_{l+1}\cdots v_{2l-3}$.
  We can check that there exists a subgraph $C_{l, l}$ in $K_2+H$, where $C_{l, l}$ is formed
   by two cycles $u'v_1v_2\cdots v_{l-1}u'$  and $u'u''v_{l}v_{l+1}\cdots v_{2l-3}u'$ which are of length $l$ and they intersect at $u'$.

   When $n_1(H)=a$ and  $n_2(H)=2l-3-a$, where $l\leqslant a  \leqslant 2l-4$, in $H$, there exist two subgraphs $P_a$ and $P_{2l-3-a}$.
  Let $P_{a}=v_1v_2\cdots v_{a}$ and $P_{2l-3-a}=v_{a+1}v_{a+2}\cdots v_{2l-3}$.
  Then $K_2+H$ contains  $C_{l, l}$ as a subgraph, where $C_{l, l}$ is formed by two cycles $u'v_1v_2\cdots v_{l-1}u'$ and $u'v_{l}\cdots v_{a}u''v_{a+1}\cdots v_{2l-3}u'$ which are of length $l$ and intersect at $u'$.

   When $n_1(H)+n_2(H)\geqslant 2l-2$, we get that  $K_2+H$ contains a subgraph $C_{l, l}$ because  there exists a subgraph $C_{l, l}$ among $K_2+H$ when $n_1(H)+n_2(H)= 2l-3$.

  \textbf{Case (2)}. $n_1(H)+n_2(H)\leqslant 2l-4$, $n_1(H)\geqslant l-1$ and $n_2(H)+n_3(H)\geqslant l-2$.

   Since $n_2(H)\geqslant n_3(H)$, we have $\lceil \frac{l}{2} \rceil-1 \leqslant n_2(H)\leqslant  l-3$.
   In $H$, there exist three subgraphs $P_{l-1}$, $P_{a}$ and $P_{l-2-a}$, where $\lceil \frac{l}{2} \rceil-1\leqslant  a\leqslant l-3$.
   We suppose $P_{l-1}=v_1v_2\cdots v_{l-1}$, $P_{a}=v_{l}v_{l+1}\cdots v_{l+a-1}$ and  $P_{l-2-a}=v_{l+a}\cdots v_{2l-3}$.
    Then $K_2+H$ contains a subgraph $C_{l, l}$, where $C_{l, l}$ is formed by two cycles $u'v_1v_2\cdots v_{l-1}u'$ and $u'v_{l}\cdots v_{l+a-1}u''v_{l+a}\cdots v_{2l-3}u'$ which are of length $l$ and intersect at $u'$.

  (2). The proof of necessity of  Claim \ref{Cl4}.

   We assume that $K_2+H$ contains a subgraph $C_{l, l}$. Let $C^1$ and $C^2$  be the two cycles of $C_{l, l}$ which are of length $l$ and intersect at a vertex.
   Three cases are considered.

 \textbf{Case (1)}. $V(C^1)\cap V(C^2)=\{u'\}$.

 Since $C^1$ and $C^2$ insect at $u'$,  we  get $u''\notin V(C^1)\cap V(C^2)$.

 \textbf{Subcase (1.1)}.  $u''\notin V(C^1)\cup V(C^2)$.

   We suppose $C^1=u'v_1v_2\cdots v_{l-1}u'$ and  $C^2=u'v_lv_{l+1}\cdots v_{2l-2}u'$.
   Then $H$ contains a subgraph $2P_{l-1}=v_1v_2\cdots v_{l-1}\cup v_lv_{l+1}\cdots v_{2l-2}$ or a subgraph $P_{2l-2}=v_1v_2\cdots v_{l-1}v_l\cdots v_{2l-2}$.
   We get $n_1(H)+n_2(H)\geqslant \min\{2l-2, 2l-1\}=2l-2$.

\textbf{Subcase (1.2)}.  $u''\notin V(C^1)$, $u''\in V(C^2)$.

  If $u'u''\in E(C^2)$,  let $C^1=u'v_1v_2\cdots v_{l-1}u'$ and  $C^2=u'u''v_lv_{l+1}\cdots v_{2l-3}u'$.
 Then $H$ contains a subgraph $P_{l-1}\cup P_{l-2}=v_1v_2\cdots v_{l-1}\cup v_lv_{l+1}\cdots v_{2l-3}$ or a subgraph  $P_{2l-3}=v_1v_2\cdots v_{l-1}v_l\cdots v_{2l-3}$.
  Thus, $n_1(H)+n_2(H)\geqslant \min\{2l-3, 2l-2\}=2l-3$.

  If  $u'u''\notin E(C^2)$, let
  $C^1=u'v_1v_2\cdots v_{l-1}u'$ and $C^2=u'v_lv_{l+1}\cdots v_{l+a}u''v_{l+a+1}\cdots v_{2l-3}\\u'$, where $0\leqslant a \leqslant l-4$.
  If $v_{l-1}v_l, v_{l+a}v_{l+a+1}\in E(H)$, then $H$ contains a subgraph $P_{2l-3}=v_1v_2\cdots v_{2l-3}$.
  Thus, $n_1(H)+n_2(H)\geqslant 2l-2$.
  If $v_{l-1}v_l\in E(H)$ and $v_{l+a}v_{l+a+1}\notin E(H)$,
  then $H$ contains a subgraph $P_{l+a}\cup P_{l-a-3}=v_1v_2\cdots v_{l-1}v_{l}\cdots v_{l+a}\cup v_{l+a+1}\cdots v_{2l-3}$, where $0\leqslant a\leqslant l-4$.
  Thus, $n_1(H)+n_2(H)\geqslant 2l-3$.
  If  $v_{l-1}v_l\notin E(H)$ and $v_{l+a}v_{l+a+1}\in E(H)$, then $H$ contains a subgraph
  $P_{l-1}\cup P_{l-2}=v_1v_2\cdots v_{l-1}\cup v_{l}v_{l+1}v_{l+2}\cdots v_{2l-3}$.
  Thus, $n_1(H)+n_2(H)\geqslant 2l-3$.
 If  $v_{l-1}v_l, v_{l+a}v_{l+a+1}\notin E(H)$, then $H$ contains a subgraph $P_{l-1} \cup P_{a+1} \cup P_{l-a-3}=v_1v_2\cdots v_{l-1}\cup v_{l}\cdots v_{l+a}\cup v_{l+a+1}\cdots v_{2l-3}$, where $0\leqslant a \leqslant l-4$.
  Thus, $n_1(H)+n_2(H)+n_3(H)\geqslant 2l-3$, $n_1(H)\geqslant l-1$ and $n_2(H)+n_3(H)\geqslant l-2$.
  When $n_1(H)+n_2(H)\geqslant 2l-3$, since $n_1(H)\geqslant n_2(H)\geqslant n_3(H)$, we get $n_1(H)\geqslant l-1$ and $n_1(H)+n_2(H)+n_3(H)\geqslant 2l-3$.
  If $n_1(H)> l-1$, since there exists a path $P_{l-1}$ in $H$, we get $n_1(H)\geqslant n_2(H)\geqslant l-1$. Thus, $n_2(H)+n_3(H)\geqslant l-1+1>l-2$.
  If $n_1(H)=l-1$, since $n_1(H)+n_2(H)\geqslant 2l-3$, we have $n_2(H)=l-2, l-1$.
   Thus, $n_2(H)+n_3(H)\geqslant l-2+1>l-2$.
 In conclusion, we get $n_1(H)+n_2(H)\geqslant 2l-3$, or $n_1(H)+n_2(H)\leqslant 2l-4$, $n_1(H)\geqslant l-1$ and $n_2(H)+n_3(H)\geqslant l-2$.

\textbf{Subcase (1.3)}.  $u''\in V(C^1)$,  $u''\notin V(C^2)$.

  By using the similar proofs as those for Subcase (1.2),  we can get the same conclusion as shown in Subcase (1.2).

\textbf{Case (2)}. $V(C^1)\cap V(C^2)=\{u''\}$.

 Since $u'$ and $u''$ are symmetric,  by using the similar proofs as those for Case (1),  we also get the same conclusion as shown in Case (1).

\textbf{Case (3)}. $V(C^1)\cap V(C^2)=\{v_i\}$, where $v_i\in V(H)$ with $1\leqslant i \leqslant n-2$.

  Since $H$ does not contain cycles and  $K_2+H$ contains $C_{l, l}$, we get that $u'$ and $u''$ are contained in $C^1$ and $C^2$, respectively.
  	Namely $u'\in V(C^1)$ and $u''\in V(C^2)$.
  Without loss of generality, let $V(C^1)\cap V(C^2)=\{v_{l-1}\}$, $C^1=u'v_1v_2\cdots v_{l-1}u'$ and $C^2=u''v_{l-1}v_l\cdots v_{2l-3}u''$.
  Then $H$ contains a subgraph $P_{2l-3}=v_1v_2\cdots v_{2l-3}$.
  Thus, $n_1(H)\geqslant 2l-3$. Namely $n_1(H)+n_2(H)\geqslant 2l-2$.

 By combining the proofs of Cases (1)--(3), we get $n_1(H)+n_2(H)\geqslant \min\{2l-3, 2l-2\}=2l-3$ or $n_1(H)+n_2(H)\leqslant 2l-4$, $n_1(H)\geqslant l-1$ and $n_2(H)+n_3(H)\geqslant l-2$. Thus, we have Claim \ref{Cl4}. \QED

   By Claim \ref{Cl1}, $ \widetilde{G}[R]$ is a union of $t$ vertex--disjoint paths, where $t\geqslant 2$.
  Let $n_i$  be the order of the $i$-th longest path of $ \widetilde{G}[R]$, where $i\in \{1, \cdots, t\}$ with $t\geqslant 2$.
     Since $n\geqslant \max\{2.16\times 10^{17}, 9\times 2^{l-1}+3, \frac{625}{32}\lfloor\frac{l-3}{2}\rfloor^2+4\}$ and
   $ n_1\geqslant n_2 \geqslant n_3\geqslant n_4 $, by Claim \ref{Cl4}, we get $n>>2(2l-4)\geqslant 2(n_1+n_2)\geqslant n_1+n_2+n_3+n_4$. Here $a>> b$  means that $a$ is much greater than $b$.
    Thus, $ \widetilde{G}[R]$ is a union of $t$ vertex--disjoint paths with $t>4$.
 \QED

\begin{claim}\label{Cl5}
	$n_1+n_2=2l-4$ and $n_1\leqslant l-2$ or $n_1+n_2=2l-4$ and  $n_2+n_3\leqslant l-3$.
\end{claim}

\textbf{Proof.} It follows from Claim \ref{Cl4} that $n_1+n_2\leqslant 2l-4$ and $n_1\leqslant l-2$ or $n_1+n_2\leqslant 2l-4$ and $n_2+n_3\leqslant l-3$.
 We prove Claim \ref{Cl5} by contradiction.
 Two cases are considered as follows.

\textbf{Case (1)}.  $n_1\leqslant l-2$.


  We assume that $n_1+n_2\leqslant 2l-6$.

  If $n_1>n_2$, let $H_1^{'}$ be the graph obtained from $ \widetilde{G}[R]$ by the $(n_2, n_t)$-transformation, 	where $t>4$.
 Then $n_1(H_1^{'})=\max\{n_1, n_2+1\}=n_1$ and $n_2(H_1^{'})=n_2+1$.
 We have $n_1(H_1^{'})+n_2(H_1^{'})=n_1+n_2+1\leqslant 2l-5$ and $n_1(H_1^{'})=n_1\leqslant l-2$.
 By  Claim \ref{Cl4},  $K_2+H_1^{'}$ does not contain $C_{l, l}$.
 By Claims \ref{Cl1} and \ref{Cl44}, we obtain $\rho(K_2+H_1^{'})>\rho(\widetilde{G})$. This contradicts the definition of $\widetilde{G}$.
 Therefore, we get  $  n_1+n_2=2l-5, 2l-4$.

 If $n_1=n_2$, since $n_1+n_2\leqslant 2l-6$, we have $n_1=n_2\leqslant l-3$.
 Let $H_2^{'}$ be the graph obtained from $ \widetilde{G}[R]$ by  the $(n_1, n_t)$-transformation, 	where $t>4$.
 We have  $n_1(H_2^{'})=n_1+1\leqslant l-2$ and  $n_2(H_2^{'})=n_2$.
 Then $n_1(H_2^{'})+n_2(H_2^{'})=n_1+n_2+1\leqslant 2l-5$ and $n_1(H_2^{'})\leqslant l-2$.
  It follows from Claim \ref{Cl4}  that  $K_2+H_2^{'}$ does not contain $C_{l, l}$.
  By Claims \ref{Cl1} and  \ref{Cl44}, we get $\rho(K_2+H_2^{'})>\rho(\widetilde{G})$. This contradicts the definition of $\widetilde{G}$.
  Thus, we have  $  n_1+n_2=2l-5, 2l-4$.

  Next, we assume that $n_1+n_2= 2l-5$.
   Since $n_1\geqslant n_2$, we have $n_1= l-2$ and  $n_2= l-3$.
   Let $H_3^{'}$ be the graph obtained from $ \widetilde{G}[R]$ by the $(n_2, n_t)$-transformation, where $t>4$.
   Then $n_1(H_3^{'})=\max\{n_1, n_2+1\}=l-2$ and  $n_2(H_3^{'})=n_2+1=l-2$.
   We get $n_1(H_3^{'})+n_2(H_3^{'})=2l-4$ and $n_1(H_3^{'})=n_1=l-2$.
   It follows from  Claim \ref{Cl4}  that $K_2+H_3^{'}$ does not contain $C_{l, l}$.
   By Claims \ref{Cl1} and \ref{Cl44}, we have $\rho(K_2+H_3^{'})>\rho(\widetilde{G})$.
   This contradicts the definition of $\widetilde{G}$. Thus, we obtain $n_1+n_2=2l-4$ and $n_1\leqslant l-2$.

\textbf{Case (2)}.  $n_2+n_3\leqslant l-3$.



 We assume that  $n_1+n_2\leqslant 2l-6$.
 Let  $H_4^{'}$ be the graph obtained from $ \widetilde{G}[R]$ by the $(n_1, n_t)$-transformation, 	where  $t>4$.
  Then $n_1(H_4^{'})=n_1+1\leqslant 2l-5-n_2$,  $n_2(H_4^{'})=n_2$ and $n_3(H_4^{'})=n_3$.
 We have  $n_1(H_4^{'})+n_2(H_4^{'})\leqslant 2l-5$ and $n_2(H_4^{'})+n_3(H_4^{'})=n_2+n_3\leqslant l-3$.
 It follows from Claim \ref{Cl4}  that  $K_2+H_4^{'}$ does not contain $C_{l, l}$.
 By Claims \ref{Cl1} and \ref{Cl44}, we have $\rho(K_2+H_4^{'})>\rho(\widetilde{G})$.
 This contradicts the definition of $\widetilde{G}$.  Therefore,  $  n_1+n_2=2l-5, 2l-4$.

  We assume that $n_1+n_2= 2l-5$.
  Let  $H_5^{'}$ be the graph obtained from $ \widetilde{G}[R]$  by  the $(n_1, n_t)$-transformation, 	where  $t>4$.
  Then $n_1(H_5^{'})=n_1+1$, $n_2(H_5^{'})=n_2$ and $n_3(H_5^{'})=n_3$.
 We have $n_1(H_5^{'})+n_2(H_5^{'})=n_1+n_2+1=2l-4$ and $n_2(H_5^{'})+n_3(H_5^{'})=n_2+n_3\leqslant l-3$.
 It follows from Claim \ref{Cl4}  that  $K_2+H_5^{'}$ does not contain $C_{l, l}$.
 By Claims \ref{Cl1} and \ref{Cl44},  we obtain  $\rho (K_2+H_5^{'})>\rho(\widetilde{G})$.  This contradicts the definition of $\widetilde{G}$.
  Thus, we obtain  $  n_1+n_2=2l-4$ and  $n_2+n_3\leqslant l-3$.

  By combining the proofs of Cases (1) and (2), we get Claim \ref{Cl5}.
  \QED


\begin{claim}\label{Cl6}
	For any $i\in\{3, \cdots, t-1 \}$, we have $n_i=n_2$, where $t>4$.
\end{claim}

\textbf{Proof.}  We suppose that there exists an $i\in\{3, \cdots, t-1 \}$  such that $n_i\leqslant n_2-1$, where $t>4$.
 Let $i_0=\min\{i~|~3\leqslant i \leqslant t-1, ~ t>4, ~n_i\leqslant n_2-1\}$.
   Two cases are considered as follows.

 \textbf{Case (1)}.  $i_0> 3$.

   Let $H_6^{'}$ be the graph obtained from $ \widetilde{G}[R]$ by  the $(n_{i_0}, n_t)$-transformation, 	where $t>4$.
   Then $n_1(H_6^{'})=n_1$, $n_2(H_6^{'})=\max\{n_2, n_{i_0}+1\}=n_2$ and $n_3(H_6^{'})=n_3$.
   It follows from Claim \ref{Cl5}  that  $n_1+n_2=2l-4$ and $n_1\leqslant l-2$ or $n_1+n_2=2l-4$ and $n_2+n_3\leqslant l-3$.
    Therefore,  $n_1(H_6^{'})+n_2(H_6^{'})=2l-4$ and $n_1(H_6^{'})=n_1\leqslant l-2$ or $n_1(H_6^{'})+n_2(H_6^{'})=2l-4$ and $n_2(H_6^{'})+n_3(H_6^{'})=n_2+n_3\leqslant l-3$.
  It follows from Claim \ref{Cl4}  that $K_2+H_6^{'}$ does not contain $C_{l, l}$.
  By Claims \ref{Cl1} and \ref{Cl44}, we get $\rho (K_2+H_6^{'})>\rho(\widetilde{G})$. This contradicts the definition of $\widetilde{G}$.

\textbf{Case (2)}.  $i_0= 3$.

 We suppose that there exists a $j$ such that $n_j<n_3\leqslant n_2-1$, where $4\leqslant j \leqslant t-1$ and $t>4$.
  Let  $H_7^{'}$  be the graph obtained from  $ \widetilde{G}[R]$ by  the $(n_j, n_t)$-transformation, 	where $4\leqslant j \leqslant t-1$.
  Then $n_1(H_7^{'})=n_1$, $n_2(H_7^{'})=n_2$ and $n_3(H_7^{'})=\max\{n_3, n_j+1\}=n_3$.
  Using the methods similar to those for  $i_0> 3$ in Case (1), we can  also get a contradiction.

  We suppose that $n_j=n_3$ for any $j\in \{4,\cdots,t-1\}$, where $t>4$.
  When $n_1+n_2=2l-4$ and $n_1\leqslant l-2$,
  let  $H_8^{'}$  be the graph obtained from  $ \widetilde{G}[R]$ by  the $(n_3, n_t)$-transformation.
  Then $n_1(H_8^{'})=n_1$, $n_2(H_8^{'})=\max\{n_2,n_3+1\}=n_2$.
  Thus we have $n_1(H_8^{'})+n_2(H_8^{'})=2l-4$ and $n_1(H_8^{'})=n_1\leqslant l-2$.
  When $n_1+n_2=2l-4$ and $n_2+n_3\leqslant l-3$, let  $H_9^{'}$  be the graph obtained from  $ \widetilde{G}[R]$ by  the $(n_1, n_2)$-transformation.
  Then $n_1(H_9^{'})=n_1+1$, $n_2(H_9^{'})=\max\{n_2-1,n_3\}=n_2-1$ and $n_3(H_9^{'})=n_3$.
  Thus we get $n_1(H_9^{'})+n_2(H_9^{'})=n_1+1+n_2-1=2l-4$ and $n_2(H_9^{'})+n_3(H_9^{'})=n_2+n_3-1\leqslant l-4<l-3$. It follows from Claim \ref{Cl4}  that neither $K_2+H_8^{'}$ nor $K_2+H_9^{'}$ contain $C_{l, l}$.
  By Claims \ref{Cl1} and \ref{Cl44}, we get $\rho (K_2+H_8^{'})>\rho(\widetilde{G})$ and $\rho (K_2+H_9^{'})>\rho(\widetilde{G})$.
  This contradicts the definition of $\widetilde{G}$.

  By combining the proofs of Cases (1) and (2) in Claim \ref{Cl6}, we get Claim \ref{Cl6}.
  \QED

 \begin{claim}\label{Cl77}
	$n_2=l-2$.
\end{claim}

\textbf{Proof.}  We suppose $n_2\neq l-2$.
  By Claim \ref{Cl5},   we have $n_1+n_2=2l-4$.
  Since $n_1\geqslant  n_2$, $n_1> l-2$.

  If $l=4$, then by Claim \ref{Cl5},  we obtain  $n_1+n_2=4$ and  $n_2+n_3\leqslant 1$.
  By Claim \ref{Cl6}, we have $n_2= n_3$.
  Thus  $n_1=4$  and $n_2=n_3=0$.
   Therefore, $\widetilde{G}\cong K_2+P_4$.
   This contradicts the condition  $n\geqslant \max\{2.16\times 10^{17}, 9\times 2^{l-1}+3, \frac{625}{32}\lfloor\frac{l-3}{2}\rfloor^2+4\}$  in Theorem \ref{MR2}.
   Thus $n_2=l-2=2$ when $l=4$.

  Next,  let $l\geqslant 5$.
  For any $i\in\{3, \cdots, t-1 \}$ with $t>4$, by Claim \ref{Cl6}, we have $n_i=n_2$.
  Since $n_2+n_3\leqslant l-3$ (by Claim \ref{Cl5}),
  we get $n_2\leqslant \lfloor\frac{l-3}{2}\rfloor$.
  Thus $n_1=2l-4-n_2 \geqslant 2l-4-\lfloor\frac{l-3}{2}\rfloor=l-2+\lceil\frac{l-1}{2}\rceil\geqslant l$.
  Since $l\geqslant 5$, $\lfloor\frac{l-3}{2}\rfloor\geqslant 1$.
  Thus $\lfloor\frac{l-3}{2}\rfloor^2\geqslant \lfloor\frac{l-3}{2}\rfloor$ and $6\lfloor\frac{l-3}{2}\rfloor^2=4\lfloor\frac{l-3}{2}\rfloor^2+2\lfloor\frac{l-3}{2}\rfloor^2\geqslant 4\lfloor\frac{l-3}{2}\rfloor+2\geqslant2l-6$. Since   $n\geqslant \frac{625}{32}\lfloor\frac{l-3}{2}\rfloor^2+4>9\lfloor\frac{l-3}{2}\rfloor^2+4$ in Theorem \ref{MR2}, we have
\begin{align}\label{q3}
	n>\left \lfloor\frac{l-3}{2}\right \rfloor^2+2\left \lfloor\frac{l-3}{2}\right \rfloor+2l-6+4\geqslant n^2_2+3n_2+n_1+2.
	\end{align}

  Since $n-2=\sum\limits _{i=1}^{t}n_i\leqslant n_1+(t-1)n_2$,  by (\ref{q3}),  we obtain $t\geqslant \frac{n-2-n_1}{n_2}+1>n_2+4$.
 By Claim \ref{Cl6}, $ \widetilde{G}[R]$ contains at least $n_2+2$ paths with $n_2$ vertices. 
 We denote $n_2+2$ paths with $n_2$ vertices by $P^1,\cdots,P^{n_2+1}, P^{n_2+2}$.
 We denote a path with $n_1$ vertices by $P_{n_1}=v_1v_2\cdots v_{n_1-1}v_{n_1}$.

  Next, we construct a graph $G^{\circ}$ as follows.
  Among  $\widetilde{G}$,  we perform  the following four  edge operations.
   (1). We delete the edge $v_{l-2}v_{l-1}$ of $P_{n_1}$ and denote the  two resulting paths by  $P_{l-2}=v_1v_2\cdots v_{l-2}$ and $P_{n_1-l+2}=v_{l-1}v_l\cdots v_{n_1}$.
   (2). We add an edge between an endpoint of $P^{n_2+2}$ and $v_{l-1}$ of $P_{n_1-l+2}$.
   (3).  We delete all the edges of $P^{n_2+1}$. Thus, we get $n_2$ isolated vertices.
   (4).  We add an edge between each isolated vertex of the $n_2$ isolated vertices and an endpoint of $P^j$ for each j $\in \{1,\cdots,n_2\}$ respectively.
    Namely, $n_2$ edges are deleted from $\widetilde{G}$ and $n_2+1$ new edges are added to $\widetilde{G}$.
   We denote the resulting  graph by $G^{\circ}$.
  Obviously,  $G^{\circ}\cong K_2+\overline{H}$, where $\overline{H}=P_{l-2}\cup P_{n_2+n_1-l+2}\cup n_2P_{n_2+1}\cup(t-n_2-4)P_{n_2}\cup P_{n_t}$.
  Since $n_1+n_2=2l-4$, we have $P_{n_2+n_1-l+2}=P_{l-2}$.
  Thus, $n_1(\overline{H})=n_2(\overline{H})=l-2$.
  Namely, $n_1(\overline{H})+n_2(\overline{H})=2l-4$ and $n_1(\overline{H})\leqslant l-2$.
  By  Claim  \ref{Cl4}, we deduce that $G^{\circ}$ does not contain $C_{l, l}$.

 For sufficiently large $n$ and any vertices $v_i, v_j\in R$ with $1\leqslant i <j \leqslant n-2$, since  $\rho(\widetilde{G})\geqslant \rho(K_{2, n-2})=\sqrt{2n-4}$, by Claim \ref{Cl3},  we get
\begin{align}\label{q4}
	\frac{4}{\rho^2(\widetilde{G})}\leqslant x_{v_i}x_{v_j}\leqslant \frac{4}{\rho^2(\widetilde{G})}+\frac{24}{\rho^3(\widetilde{G})}+\frac{36}{\rho^4(\widetilde{G})}<\frac{4}{\rho^2(\widetilde{G})}+\frac{25}{\rho^3(\widetilde{G})}.
	\end{align}
 Since $n\geqslant \frac{625}{32}\lfloor\frac{l-3}{2}\rfloor^2+4$, we obtain $n_2\leqslant \lfloor\frac{l-3}{2}\rfloor< \frac{4}{25}\sqrt{2n-4}\leqslant \frac{4}{25}\rho(\widetilde{G})$.
 By (\ref{q2}) and (\ref{q4}),
 \begin{align*}
	\rho(G^{\circ})-\rho(\widetilde{G})&\geqslant \frac{\bm{x}^\mathrm{T}(A(G^{\circ})-A(\widetilde{G}))\bm{x}}{\bm{x}^\mathrm{T}\bm{x}} \nonumber\\
	&>\frac{2}{\bm{x}^\mathrm{T}\bm{x}}\left(\frac{4(n_2+1)}{\rho^2(\widetilde{G})}-\frac{4n_2}{\rho^2(\widetilde{G})}-\frac{25n_2}{\rho^3(\widetilde{G})}\right)>0.
\end{align*}
Namely, $\rho(G^{\circ})>\rho(\widetilde{G})$.
This contradicts the definition of $\widetilde{G}$.
Thus we get Claim \ref{Cl77}.
\QED

 By combination of the proofs of Claims \ref{Cl1}--\ref{Cl77}, we have
 $n_1=2l-4-n_2=l-2$, $n_i=n_2=l-2$ with $i\in\{3, \cdots , t-1\}$ and  $n_t\leqslant l-2$.
 This means $ \widetilde{G}[R]\cong H(l-2, l-2)$.
 Thus, $\widetilde{G}\cong K_2+H(l-2, l-2)$.
 Namely,  Theorem \ref{MR2} holds. \QED

\subsection{Extremal spectral results of planar graphs without
 $\mathrm{Theta}$ graph}

   In this section,  among the set of planar graphs on $n$ vertices without any member of $\Theta_k$,
   for sufficiently large $n$, we characterize the extremal graph  to $\textrm{spex}_{P}(n, \Theta_k)$,  where $k\geqslant4 $.
  The results are shown in Theorems \ref{MR3} and \ref{MR3aa}.

\begin{theorem}\label{MR3}
   	  Let $G$ be a planar graph on $n$ vertices without $C_3\cdot C_{3}$.
    When $n\geqslant 2.16\times 10^{17}$,
      we have $\rho(G)\leqslant \rho(K_{2, n-2})$ with the equality if and only if $G\cong K_{2, n-2}$.
	\end{theorem}

   \textbf{Proof.}  We suppose that $G^{\diamond}$ is the  extremal graph  to $\textrm{spex}_{P}(n, C_3\cdot C_{3})$,   where $n\geqslant 2.16\times 10^{17}$.
      Since $C_3\cdot C_{3}\subseteq2K_1+P_{n/2}$ and $C_3\cdot C_{3} \nsubseteq K_{2, n-2}$, by Lemma \ref{PF1}, we get $K_{2, n-2}\subseteq G^{\diamond}$.
     In $K_{2, n-2}$, if we add an edge between two vertices which are not adjacent  in $K_{2, n-2}$, then the resulting planar graph should contain $C_3\cdot C_{3}$.
    Therefore, we have  $G^{\diamond}\cong K_{2, n-2}$.
 \QED

 \begin{theorem}\label{MR3aa}
   Let $G$ be a planar graph on $n$ vertices and $G$  does not contain any member of $\Theta_k$  as a subgraph, where $k\geqslant5$.
   When $n\geqslant \max\{2.16\times 10^{17}, 9\times 2^{\lfloor \frac{k-1}{2}\rfloor}+3, \frac{625}{32}\lfloor\frac{k-3}{2}\rfloor^2+2\}$,
 we have $\rho(G)\leqslant \rho(K_2+H(\lceil\frac{k-3}{2}\rceil, \lfloor\frac{k-3}{2}\rfloor))$ with the equality if and only if $G\cong K_2+H(\lceil\frac{k-3}{2}\rceil, \lfloor\frac{k-3}{2}\rfloor)$.
\end{theorem}

 \textbf{Proof.}  We assume that $\overline{G}$ is the  extremal graph  to $\textrm{spex}_{P}(n, \Theta_k)$,  where  $k\geqslant5$ and $n\geqslant \max\{2.16\times 10^{17}, 9\times 2^{\lfloor \frac{k-1}{2}\rfloor}+3, \frac{625}{32}\lfloor\frac{k-3}{2}\rfloor^2+2\}$.
    Let $G''$ be an arbitrary graph among $\Theta_k$, where $k\geqslant 5$. Since $n\geqslant 9\times 2^{\lfloor \frac{k-1}{2}\rfloor}+3> 2k=2|V(G'')|$,
    we have $G''\subseteq2K_1+P_{n/2}$ and $G'' \nsubseteq K_{2, n-2}$. By Lemma \ref{PF1}, we get $K_{2, n-2}\subseteq \overline{G}$.
    Thus, $\rho(\overline{G})\geqslant \rho(K_{2, n-2})=\sqrt{2n-4}$.
    By Perron-Frobenius theorem,  there exists a positive eigenvector $\bm{x}=(x_1,  x_2,  \cdots,  x_n)^{\mathrm{T}}$ corresponding to $\rho(\overline{G})$.
    By the proofs of Lemma \ref{PF1} in Ref.\ \cite{https://doi.org/10.1002/jgt.23084} (see Pages 6-7),  among $\bm{x}$,  there exist two components with larger values.
    In $\overline{G}$,  let $u_1$ and $u_2$ be the two vertices corresponding to the two larger components,  where
     $x_{u_1}=\max\{x_i: 1\leqslant i \leqslant n\}=1$. Let $A=\{u_1, u_2\}$ and $B=V(\overline{G})\setminus A=\{w_1, \cdots , w_{n-2}\}$.
   Since $K_{2, n-2}\subseteq \overline{G}$,  we have $B=N(u_1)\cap N(u_2)$ (by the proofs of Lemma \ref{PF1} in Ref.\ \cite{https://doi.org/10.1002/jgt.23084}).

  To obtain Theorem \ref{MR3aa}, we introduce Claims \ref{C2222l7}--\ref{Cxxl899}.

  \begin{claim}\label{C2222l7}
	$\overline{G}[B]$ is a union of $t$ vertex--disjoint paths, where $t\geqslant 2$.
\end{claim}

\textbf{Proof.}  By using the same methods as those for the proofs of Claim \ref{Cl1}, we obtain Claim \ref{C2222l7}.
\QED

 For any $j\in\{1, \cdots, t\}$ with  $t\geqslant 2$, let $n_j$ be the order of the $j$-th longest path of $\overline{G}[B]$. Namely,  $n_1\geqslant n_2\geqslant \cdots \geqslant n_t$.

\begin{claim}\label{Cl7}
	$u_1u_2\in E(\overline{G})$.
\end{claim}

\textbf{Proof.}  First, we prove $n_1\leqslant k-4$.
 Otherwise, we suppose  $n_1\geqslant k-3$.
  Then $\overline{G}[B]$ contains $P_{k-3}=w_1w_2\cdots w_{k-3}$.
  Therefore,  $\overline{G}$ contains $C_{a+1}\cdot C_{k+1-a}$ as a subgraph,
  where $C_{a+1}\cdot C_{k+1-a}$ is formed by  two cycles $u_1w_1\cdots w_au_1$ and $u_1w_a\cdots w_{k-3}u_2w_{k-2}u_1$ with $2\leqslant a \leqslant \lfloor\frac{k}{2}\rfloor$ and $k\geqslant 5$.
  This contradicts the fact that $\overline{G}$ does not contain any member of $\Theta_k$ as a subgraph.
   Thus, we get $n_1\leqslant k-4$.

   Next, we prove $n_1+n_2\leqslant k-3$.
  Otherwise, we suppose  $n_1+n_2\geqslant  k-2$.
  When $k=5$, $n_1+n_2\geqslant  3$.
 Since $n_1\leqslant 5-4=1$, we have $n_2\geqslant 2$. This contradicts the fact that $n_1\geqslant n_2$.
 When $k\geqslant 6$,  since $n_1\leqslant k-4$, $n_1+n_2\geqslant  k-2$ and  $n_1\geqslant n_2$, we have $2\leqslant n_2\leqslant n_1\leqslant k-4$.
 Thus, $\overline{G}[B]$ contains $P_a\cup P_{k-2-a}=w_1\cdots w_a\cup w_{a+1}\cdots w_{k-2}$ as a subgraph, where $2\leqslant a \leqslant \lfloor\frac{k}{2}\rfloor-1$.
  We can check that  $\overline{G}$ contains $C_{a+1}\cdot C_{k+1-a}$,  where $C_{a+1}\cdot C_{k+1-a}$ is formed by  two cycles $u_1w_1\cdots w_au_1$ and $u_1w_au_2w_{a+1}\cdots w_{k-2}u_1$ with $2\leqslant a \leqslant \lfloor\frac{k}{2}\rfloor-1$.
   This contradicts the fact that $\overline{G}$ does not contain any member of $\Theta_k$  as a subgraph.
  Thus, $n_1+n_2\leqslant k-3$.

 We suppose $u_1u_2\notin E(\overline{G})$.
 Since $\overline{G}[B]$ is a union of $t$ vertex-disjoint paths  (by Claim \ref{C2222l7}),
  there exists an  $i\in\{1, \cdots, n-2\}$ such that
 $w_iw_{i+1}\notin E(\overline{G}[B])$ (we define $w_{n-1}=w_1$), where $t\geqslant 2$.
 We can check that $u_1w_iu_2w_{i+1}u_1$ is a face of $\overline{G}$, where $1 \leqslant i \leqslant n-2$.
 Next, we construct a graph $G^{\bigtriangleup}=\overline{G}+\{u_1u_2\}$,  where $u_1u_2$ passes through the surface $u_1w_iu_2w_{i+1}u_1$ with $1 \leqslant i \leqslant n-2$.
  Obviously, $G^{\bigtriangleup}$ is a planar graph.
  We will prove that $G^{\bigtriangleup}$ does not contain any member of $\Theta_k$  as a subgraph, where $k\geqslant5$.

 We suppose that $G^{\bigtriangleup}$ contain  a member of $\Theta_k$  as a subgraph, where $k\geqslant5$.
 In $G^{\bigtriangleup}$,  there exist two cycles (denoted by $\overline{C^1}$ and  $\overline{C^2}$) which share a common edge and
  $u_1u_2\in E(\overline{C^1})\cup E(\overline{C^2})$.
 Let $|V(\overline{C^1})|=r$ and  $|V(\overline{C^2})|=s$, where $r$ and $s$ satisfy $3\leqslant r\leqslant \lfloor\frac{k}{2}\rfloor+1$,
  $\lceil \frac{k}{2}\rceil+1 \leqslant s\leqslant k-1$ and  $r+s=k+2$.
  We suppose $E(\overline{C^1})\cap E(\overline{C^2})=\{e\}$.
  Two cases are considered.

  \textbf{Case (1)}.  $e=u_1u_2$.

   We assume that $\overline{C^1}=u_1u_2w_1\cdots w_{r-2}u_1$ and  $\overline{C^2}=u_1u_2w_{r-1}\cdots w_{r+s-4}u_1$.
   Then,  in $\overline{G}[B]$, there exist two vertex--disjoint paths: $P_{r-2}\cup P_{s-2}=w_1\cdots w_{r-2} \cup w_{r-1}\cdots w_{r+s-4}$
    or a path with $r+s-4=k-2$ vertices: $P_{k-2}=w_1w_2\cdots w_{k-2}$.
   Thus, we have $n_1+n_2\geqslant \min\{r-2+s-2, k-1\} =\min\{k+2-4, k-1\}=k-2$.
   This contradicts  $n_1+n_2\leqslant k-3$.

\textbf{Case (2)}.  $e \neq u_1u_2$.

\textbf{Subcase (2.1)}. $e\notin E(\overline{G}[B])$.

 We assume $\overline{C^1}=u_1u_2w_1\cdots w_{r-2}u_1$.
 Without loss of generality, we suppose  $ e=w_{r-2}u_1$.
 Let $\overline{C^2}=u_1w_{r-2}w_{r-1}\cdots w_{r+s-4}u_1$.
  In $\overline{G}[B]$,  there exists a path with $r+s-4=k-2$ vertices: $P_{k-2}=w_1w_2\cdots w_{r+s-4}$.
 Thus, $n_1+n_2\geqslant k-1$. This contradicts $n_1+n_2\leqslant k-3$.

\textbf{Subcase (2.2)}. $e\in E(\overline{G}[B])$.

 We assume $\overline{C^1}=u_1u_2w_1\cdots w_{r-2}u_1$.
 We suppose $e= w_{r-3}w_{r-2}$.
 Let $\overline{C^2}=w_{r-3}w_{r-2}w_{r-1}\cdots w_{r+s-4}w_{r-3}$.
 Then $\overline{G}[V(\overline{C^2})\cup \{u_1, u_2\}]$ contains $K_5$-minor.
  This contradicts the fact that $\overline{G}$ is a planar graph.

  By combining the proofs of Cases (1) and (2), we obtain that $G^{\bigtriangleup}$ is a planar graph and $G^{\bigtriangleup}$ does not contain any member of  $\Theta_k$ as a subgraph, where $k\geqslant 5$.
 Since $\overline{G}\subset G^{\bigtriangleup}$, by Lemma \ref{liqiao},  we have $\rho(\overline{G}) < \rho(G^{\bigtriangleup})$.
 This contradicts the definition of $\overline{G}$. Thus, we have $u_1u_2\in E(\overline{G})$.\QED
\begin{claim}\label{cl888}
	For any $w_i\in B$, $x_{w_i}\in[\frac{2}{\rho(\overline{G})}, \frac{2}{\rho(\overline{G})}+\frac{6}{\rho^2(\overline{G})}]$,
	where $1\leqslant i \leqslant n-2$.
\end{claim}

\textbf{Proof.}
 By the methods similar to those for the proofs of Claim \ref{Cl3}, we obtain Claim \ref{cl888}.\QED

\begin{claim}\label{cl88}
	Let $H$ and $H^*$ be the two graphs as shown in
		Definition \ref{df1}, $n \geqslant max\{2.16\times 10^{17}, 9\times 2^{s_2+1}+3,2|V(F)|\}$, where $F \in \Theta_k$ and $k \geqslant 5$. When $ \overline{G}[B]\cong H$, we have $\rho(K_2+H^*)>\rho(\overline{G})$.
	\end{claim}

\textbf{Proof.} By Claims \ref{C2222l7}--\ref{cl888} and using the same methods as those for the proofs of Lemma 3.2 in Ref.\ \cite{https://doi.org/10.1002/jgt.23084} (see Pages 10-13), we obtain Claim \ref{cl88}.\QED

 \begin{claim}\label{Cl8}
	Let $H$ be a  union of $t$  vertex--disjoint paths, where $t\geqslant 2$.
     For any $i\in \{1, \cdots, t\}$, let $n_i(H)$ be the order of the $i$-th longest path of $H$.
     When $k\geqslant 5$, $K_2+H$ does not contain any member of $\Theta_k$ as a subgraph  if and only if $n_1(H)+n_2(H)\leqslant k-3$.
\end{claim}

\textbf{Proof.} Let $V(K_2)=\{u_1, u_2\}$ and  $V(H)=\{w_i: 1 \leqslant i \leqslant |V(H)|\}$.
  By the methods similar to the proofs of the second paragraph in Claim \ref{Cl7} for $n_1+n_2\leqslant k-3$, we obtain the necessity of  Claim \ref{Cl8}.
  Next, we prove the sufficiency of Claim \ref{Cl8}.
  Namely, if $n_1(H)+n_2(H)\leqslant k-3$, then we will prove that $K_2+H$ does not contain any member of $\Theta_k$ as a subgraph.
  This is equivalent to prove that if $K_2+H$ contains a member of $\Theta_k$ as a subgraph, then $n_1(H)+n_2(H)\geqslant k-2$.

   We suppose that $K_2+H$ contain a member of $\Theta_k$ as a subgraph, where $k\geqslant 5$.
   Then in $K_2+H$ there exist two cycles (denoted by $\overline{C^1}$ and $\overline{C^2}$) which share a common edge.
  Let $|V(\overline{C^1})|=r$ and  $|V(\overline{C^2})|=s$,  where $r$ and $s$ satisfy $3\leqslant r\leqslant \lfloor\frac{k}{2}\rfloor+1$, $\lceil \frac{k}{2}\rceil+1 \leqslant s\leqslant k-1$ and  $r+s=k+2$.
   Two cases are considered as follows.

\textbf{Case (1)}.  $u_1u_2 \in E(\overline{C^1})\cup E(\overline{C^2})$.

  By the methods similar to the proofs of Claim \ref{Cl7} for Cases (1) and (2), we get $n_1(H)+n_2(H)\geqslant \min\{k-2, k-1\}=k-2$.

\textbf{Case (2)}.  $u_1u_2 \notin E(\overline{C^1})\cup E(\overline{C^2})$.

 Let $E(\overline{C^1})\cap E(\overline{C^2})=\{e\}$.

 \textbf{Subcase (2.1)}. $e\in E(H)$.

   Let $\overline{C^1}=u_1w_1\cdots w_{r-1}u_1$. Without loss of generality, we suppose $ e=w_{r-2}w_{r-1}$.
  Let $\overline{C^2}=u_2w_{r-2}w_{r-1}\cdots w_{r+s-4}u_2$.
   In $H$, there exists a path with $r+s-4=k-2$ vertices: $P_{k-2}=w_1w_2\cdots w_{k-2}$.
   Thus, $n_1(H)+n_2(H)\geqslant k-1$.

 \textbf{Subcase (2.2)}. $e\notin E(H)$.

  Let $\overline{C^1}=u_1w_1\cdots w_{r-1}u_1$. We suppose $ e=u_1w_{r-1}$.
  If $\overline{C^2}=u_1w_{r-1}\cdots w_{r-1+a}u_2\\w_{r+a}\cdots w_{r+s-4}u_1$,
  then $H$ contains a path with $r+s-4=k-2$ vertices: $P_{k-2}=w_1w_2\cdots w_{k-2}$ for $w_{r-1+a}w_{r+a}\in E(H)$ or
  $H$ contains $P_{r-1+a}\cup P_{k-2-(r-1+a)}=w_1w_2\cdots w_{r-1+a} \cup w_{r+a}w_{r+a+1}\cdots w_{k-2}$ for $w_{r-1+a}w_{r+a}\notin E(H)$, where $0\leqslant a \leqslant s-4$.
  If $\overline{C^2}=u_1w_{r-1}w_{r}\cdots w_{r+s-3}u_1$, then $H$ contains a path with $r+s-3=k-1$ vertices: $P_{k-1}=w_1w_2\cdots w_{k-1}$.
  Thus, $n_1(H)+n_2(H)\geqslant \min\{k-1, k-2\}=k-2$.

  By combining the proofs of  Cases (1) and (2), we get  $n_1(H)+n_2(H)\geqslant k-2$.\QED

  It is noted that $\overline{G}[B]$ is  a union of $t$ vertex--disjoint paths, where $t\geqslant 2$.
   By Claim \ref{Cl8}, $n_1+n_2\leqslant k-3$. 
 Since $n_1\geqslant n_2$,  we get $n_2\leqslant \lfloor \frac{k-3}{2}\rfloor$.
 Since $k\geqslant 5$, then $\lfloor \frac{k-3}{2}\rfloor\geqslant 1$.
 Thus, $\lfloor \frac{k-3}{2}\rfloor^2\geqslant \lfloor \frac{k-3}{2}\rfloor$ and $3\lfloor \frac{k-3}{2}\rfloor^2\geqslant 2\lfloor \frac{k-3}{2}\rfloor+1\geqslant k-3$.
  Since  $n\geqslant \frac{625}{32}\lfloor\frac{k-3}{2}\rfloor^2+2>6\lfloor \frac{k-3}{2}\rfloor^2+2$ in Theorem \ref{MR3aa}, we get
 \begin{align}\label{q5}
	n>\left \lfloor\frac{k-3}{2}\right \rfloor^2+2\left\lfloor\frac{k-3}{2}\right\rfloor+k-3+2\geqslant n^2_2+3n_2+n_1+2.
\end{align}

 Since $n-2=\sum\limits _{i=1}^{t}n_i\leqslant n_1+(t-1)n_2$, by (\ref{q5}), we get $t\geqslant \frac{n-2-n_1}{n_2}+1>n_2+4$.

  \begin{claim}\label{Cxxl8} $n_1+n_2=k-3$.
   \end{claim}

 \textbf{Proof.}  We suppose $n_1+n_2\leqslant k-4$.
   Let $H_1^{*}$ be the graph obtained from $\overline{G}[B]$ by the $(n_1, n_t)$-transformation, where $t\geqslant n_2+4$.
   Then $n_1(H_1^{*})=n_1+1\leqslant k-3-n_2$ and  $n_2(H_1^{*})=n_2$.
    Thus $n_1(H_1^{*})+n_2(H_1^{*})\leqslant k-3$.
    By Claim \ref{Cl8},  we get that $K_2+H_1^{*}$ does not contain any member of $\Theta_k$ as a subgraph.
    By Claims \ref{C2222l7} and \ref{cl88}, $\rho(K_2+H_1^{*})>\rho(\overline{G})$.
   This contradicts the definition of $\overline{G}$. Thus, we get $n_1+n_2\geqslant n-3$.
   Furthermore, since  $n_1+n_2\leqslant  n-3$ (by Claim \ref{Cl8}), 
    we get Claim \ref{Cxxl8}. \QED

   \begin{claim}\label{Cddddl8}
   We have $n_i=n_2$ for any $i\in\{3, \cdots , t-1\}$, where $t\geqslant n_2+4$.
   \end{claim}

   \textbf{Proof.}  By the methods similar to those for the proofs of Claim \ref{Cl6}, we get   Claim \ref{Cddddl8}.
  \QED

   Since $t\geqslant n_2+4$, $\overline{G}[B]$ contains at least $n_2+2$ paths with $n_2$ vertices. We denote $n_2+2$ paths with $n_2$ vertices by $P^1,\cdots,P^{n_2+1}, P^{n_2+2}$. 
   Let $P_{n_1}=w_1w_2\cdots w_{n_1-1}w_{n_1}$ be a path of order $n_1$.

  \begin{claim}\label{Cxxl899} $n_2=\left \lfloor\frac{k-3}{2}\right \rfloor$.
   \end{claim}

   \textbf{Proof.}   We suppose $n_2\leqslant \left \lfloor\frac{k-5}{2}\right \rfloor$.
  Since $n_1+n_2=k-3$, we get $n_1\geqslant k-3-\left \lfloor\frac{k-5}{2} \right\rfloor=\left\lceil \frac{k-1}{2} \right\rceil$.
  We construct a graph $G^{\bigtriangledown}$ as follows.
   Among  $\overline{G}$,  we perform  the following three  edge operations.
 (1).  We delete the edge $w_{n_1-1}w_{n_1}$ from $P_{n_1}$,  and add an edge between an endpoint of $P^{n_2+2}$ and  $w_{n_1}$.
 (2). We delete all the edges of $P^{n_2+1}$. Then we get $n_2$ isolated vertices.
  (3). We add an edge between each isolated vertex of the $n_2$ isolated vertices and an endpoint of $P^j$ for each j $\in \{1,\cdots,n_2\}$, respectively.
    Namely, $n_2$ edges are deleted from $\overline{G}$ and $n_2+1$ new edges are added to $\overline{G}$.
    We denote the resulting graph by $G^{\bigtriangledown}$.
  By the construction of $G^{\bigtriangledown}$, $G^{\bigtriangledown}\cong K_2+H^{*}$,  where $H^{*}=P_{n_1-1}\cup(n_2+1)P_{n_2+1}\cup(t-n_2-4)P_{n_2}\cup P_{n_t}$.
  We can check that $n_1(H^{*})=\max\{n_1-1, n_2+1\}=n_1-1$ and  $n_2(H^{*})=n_2+1$. Thus $n_1(H^{*})+n_2(H^{*})=k-3$.
  Therefore, by Claim \ref{Cl8}, we get that  $G^{\bigtriangledown}$ does not contain any member of $\Theta_k$ as a subgraph, where $k\geqslant 5$.

 Since  $n $ is sufficiently large  in Theorem \ref{MR3aa} and $\rho(\overline{G})\geqslant \rho(K_{2, n-2})=\sqrt{2n-4}$, for any vertices $w_i,  w_j\in B$ with  $1\leqslant i<j\leqslant n-2$, by Claim \ref{cl888}, we have
 \begin{align}\label{q6}
	\frac{4}{\rho^2(\overline{G})}\leqslant x_{w_i}x_{w_j}\leqslant \frac{4}{\rho^2(\overline{G})}+\frac{24}{\rho^3(\overline{G})}+\frac{36}{\rho^4(\overline{G})}<\frac{4}{\rho^2(\overline{G})}+\frac{25}{\rho^3(\overline{G})}.
\end{align}

 Since  $n\geqslant \frac{625}{32}\lfloor\frac{k-3}{2}\rfloor^2+2$ in Theorem \ref{MR3aa},  we have  $n_2\leqslant \lfloor\frac{k-5}{2}\rfloor< \lfloor\frac{k-3}{2}\rfloor\leqslant \frac{4}{25}\sqrt{2n-4}\leqslant \frac{4}{25}\rho(\overline{G})$.
 By (\ref{q2}) and (\ref{q6}), we get
\begin{align*}
	\rho(G^{\bigtriangledown})-\rho(\overline{G})&\geqslant \frac{\bm{x}^\mathrm{T}(A(G^{\bigtriangledown})-A(\overline{G}))\bm{x}}{\bm{x}^\mathrm{T}\bm{x}} \nonumber\\
	&>\frac{2}{\bm{x}^\mathrm{T}\bm{x}}\left(\frac{4(n_2+1)}{\rho^2(\overline{G})}-\frac{4n_2}{\rho^2(\overline{G})}-\frac{25n_2}{\rho^3(\overline{G})}\right)>0.
\end{align*}
Namely, $\rho(G^{\bigtriangledown})>\rho(\overline{G})$.
This contradicts the definition of $\overline{G}$.
Thus, $n_2=\left \lfloor\frac{k-3}{2}\right \rfloor$.
  \QED

 By Claims \ref{Cxxl8} and \ref{Cxxl899}, we have  $n_1=k-3-n_2=\left \lceil\frac{k-3}{2}\right \rceil$.
 It follows from Claim \ref{Cddddl8} that  $n_i=n_2=\left \lfloor\frac{k-3}{2}\right \rfloor$ with $i\in\{3, \cdots , t-1\}$ and $n_t\leqslant \left \lfloor\frac{k-3}{2}\right \rfloor$.
This means that $\overline{G}[B]\cong H(\lceil\frac{k-3}{2}\rceil, \lfloor\frac{k-3}{2}\rfloor)$.
Therefore, we get $\overline{G}\cong K_2+H(\lceil\frac{k-3}{2}\rceil, \lfloor\frac{k-3}{2}\rfloor)$.
Thus, we obtain Theorem \ref{MR3aa}. \QED

%

\begin{acknowledgments}
 The work was supported by the Natural Science Foundation of Shanghai under the grant number 21ZR1423500.
\end{acknowledgments}

\bibliographystyle{elsevier_citation_order}  
\bibliography{bibliography20240129}
\end{document}